\documentclass[11pt]{amsart}
\usepackage[utf8]{inputenc}   
\usepackage[T1]{fontenc}
\usepackage{amsmath}     
\usepackage{amssymb}
\usepackage{mathrsfs}
\usepackage{dsfont}        
\usepackage{bbm}

\usepackage{amsthm}
\usepackage[ngerman, english]{babel}            
\usepackage{graphicx}
\usepackage[english] {babel}
\usepackage[a4paper,margin=3cm]{geometry}
\usepackage{graphicx}
\usepackage{tikz}

\usepackage[all]{xy}  

\usepackage{color}

\usepackage[toc,page,title]{appendix}  
\usepackage{hyperref}
\hypersetup{
    colorlinks=true,
    linkcolor=blue,
    filecolor=magenta,      
    urlcolor=cyan,
  }

\DeclareMathOperator{\id}{id}

\DeclareMathOperator{\tr}{tr}

\DeclareMathOperator{\m}{\mathfrak{m}}

\DeclareMathOperator{\R}{\mathbb{R}}
\DeclareMathOperator{\N}{\mathbb{N}}

\DeclareMathOperator{\Z}{\mathbb{Z}}

\DeclareMathOperator{\HH}{\mathcal{H}}
\DeclareMathOperator{\LL}{\mathcal{L}}
\DeclareMathOperator{\WW}{\mathcal{W}}
\DeclareMathOperator{\VV}{\mathcal{V}}
\DeclareMathOperator{\UU}{\mathcal{U}}
\DeclareMathOperator{\BB}{\mathcal{B}}
\DeclareMathOperator{\FF}{\mathcal{F}}
\DeclareMathOperator{\AAA}{\mathcal{A}}
\DeclareMathOperator{\TT}{\mathcal{T}}

\DeclareMathOperator{\A}{\mathring{A}}

\DeclareMathOperator{\Ric}{Ric}
\DeclareMathOperator{\Sc}{Sc}
\DeclareMathOperator{\Rm}{Rm}

\DeclareMathOperator{\Div}{div}
\DeclareMathOperator{\diam}{diam}
\renewcommand{\d}{\mathrm{d}}

\DeclareMathOperator{\Hess}{Hess}
\DeclareMathOperator{\Vol}{Vol}
\DeclareMathOperator{\dist}{dist}
\DeclareMathOperator{\inj}{inj}

\newcommand{\suchthat}{\, \middle| \, }

\theoremstyle{plain}
\newtheorem{theorem}{Theorem}[section]
\newtheorem{lemma}[theorem]{Lemma}
\newtheorem{prop}[theorem]{Proposition}
\newtheorem{corollary}[theorem]{Corollary}

\theoremstyle{definition}
\newtheorem{definition}[theorem]{Definition}

\theoremstyle{remark}
\newtheorem*{remark}{Remark}

\numberwithin{equation}{section}  

\title{Concentration of Small Hawking Type Surfaces}
\author{Alexander Friedrich } 
\date{September 5, 2019}
\address{University of Potsdam, Institute of Mathematics, Karl-Liebknecht-Str. 24-25,
D-14476 Potsdam OT Golm}
\thanks{The author was supported by the DFG project ME3816/1-2. Further, the author would like to thank Jan Metzger for his guidance and patience during the authors Ph.D. from which this article developed.}

\begin{document}

\begin{abstract}
We investigate the Hawking energy of small surfaces in space times without symmetry assumptions by introducing the notion of Hawking type functionals. In particular, we find that Hawking type functionals are generalized Willmore functionals which allows us to find area constrained, minimizing, immersed, haunted bubble trees. These bubble trees are smooth spheres provided their area is small enough.  

Following a similar analysis of the Willmore functional conducted by T. Lamm and J. Metzger we characterize the concentration points of area constrained, critical surfaces for Hawking type functionals and the Hawking energy. Moreover, we determine their expansion on small surfaces.
\end{abstract}

\thispagestyle{empty}

\maketitle

\pagestyle{headings}

\section{Introduction}
In this paper we investigate the Hawking energy in the context of generalized Willmore functionals as introduced in \cite{Me_minimizers_gen}. 

The Hawking energy is a quasi local energy functional used in general relativity.
It was first proposed by S. W. Hawking in \cite{Hawking_energy} as a measure for the classical energy as well as the gravitational energy of an isolated system. Under the premise that energy determines the geometry of space time, the idea is to measure the bending of light rays across a spherical surface and to compare it to the flat case. For more on the Hawking energy and other quasi local energies see \cite{Szabados_quasi_local_review}.

General Relativity  is modeled on a Lorentz manifold, the most important aspects of which we briefly summarize.
Let $(N^4,h)$ be a four dimensional Lorentz manifold and let $(M_t^3,g_t)$ be an oriented, space like foliation of $N$. That is for every $t \in \R$, $(M_t,g_t)$ is a Riemannian manifold, where $g_t$ is the restriction of $h$ to $M_t$, which we interpret as equal time slice. We will focus on a given leaf and thus drop the $t$ dependence. The second fundamental form $K$ of $M$ in $N$ is given by
\[K(X,Y) := h(\nabla^N_X n , Y), \] 
where $X$ and $Y$ are vector fields of $M$ and $n$ is the (time like) normal vector of $M$.\\
Further, consider an immersed Riemann surface
$\phi: S \to \Sigma \subset M$ with induced metric $\gamma = \phi^* g$. Its area is denoted by $|\Sigma|_\gamma = \int_\Sigma \, \d \mu_\gamma$. If the context is clear we will drop the metric dependency from the notation.\\
The second fundamental form of $\Sigma$ in $M$ defined by 
\[A(X,Y):= g(\nabla^M_X \nu, Y), \]
where $X$ and $Y$ are tangent vector fields of $\Sigma$, $\nu$ is the normal vector field of $\Sigma$ in $M$ and $\nabla^M$ is the Levi-Civita connection on $M$. It is useful to decompose $A$ into its trace free part $\A$ and its  trace   $A := \A + \frac{1}{2} \gamma H $, where $H:= \tr_\Sigma A$ is the mean curvature of $\Sigma$.

Analogous to  the mean curvature of $\Sigma$ in $M$, we define the mean curvature of $\Sigma$ with respect to $K$ to be
\begin{align*}
P&:= \tr_\Sigma K = \tr_M K - K(\nu,\nu).
\end{align*}
Then the mean curvature vector of $\Sigma$ in $N$ is given by $\vec{H} := H \nu + P n$.

The Hawking energy of $\Sigma$ is defined as
\begin{align*}
 \mathcal{E}[\Sigma] &= \sqrt{\frac{|\Sigma|}{16 \pi}} \left(1-   \frac{1}{16 \pi} \int_\Sigma |\vec{H}|_h^2 \, \d\mu \right) \\
 &= \sqrt{\frac{|\Sigma|}{16 \pi}} \left(1-   \frac{1}{16 \pi} \int_\Sigma H^2 - P^2 \, \d\mu \right) .
\end{align*} 
Clearly, minimizing the functional $\int_\Sigma |\vec{H}|_h^2 \, \d\mu$ under area constraint amounts to  maximizing the Hawking energy under area constraint. 
Here we take a more general approach and investigate \emph{Hawking type} functionals of the following form. Let $L: TM \to \R$ be given and define
\[ \HH[\Sigma] := \HH_L[\Sigma] := \WW[\Sigma] + \int_\Sigma L(x,\nu) \, \d\mu, \]
 where  $\WW[\Sigma]$ is the the Willmore functional.
\begin{align*}
\WW[\Sigma] &:= \frac{1}{4} \int_\Sigma H^2 d\mu 
\end{align*}
In Section \ref{section_Hawking_intro} we establish the existence and regularity of area constrained minimizers of Hawking type functionals and therefore area constrained maximizers of $\mathcal{E}$. 
\begin{theorem}
Let $(M,g)$ be $C_B$-bounded and let  $\HH_L$ be a Hawking type functional for a smooth and bounded $L$. Then $\HH_L$ is a generalized Willmore functional in the sense Definition \ref{def_gen_will}. Moreover,
\begin{enumerate}
	\item if $M$ is compact, then the inifimum of $\HH_L$  among haunted, branched, immersed bubble trees with area $a$ is attained  for any $a>0$.	
	 Additionally, any area constrained critical point of $\HH_L$ is smooth away from finitely many points.
\item
There is a constant $a_0(L, C_B)>0$ such that all $\Sigma_a$ realizing the infimum of $\HH_L$ with area $a\leq a_0$ are embedded spheres, contained in a normal coordinate neighborhood and satisfy
\[| \HH[\Sigma_a] - 4\pi| \leq C(L, C_B) a. \]
\end{enumerate}
\end{theorem}

In Section \ref{section_concentration} we investigate Hawking type functionals in the spirit of T. Lamm and J. Metzger in \cite{LMI}. That is, we calculate expansions on small spheres and characterize concentration points, i.e. points  in the ambient manifold around which there exit critical, area constrained, spherical surfaces $\Sigma_r$ in any neighborhood $B_r(p)$. In particular, we prove the following results for the Hawking energy.
\begin{theorem}
\label{thm_into_concentration}
Let $(M,g)$ be $C_B$-bounded and let $\HH[\Sigma] = \frac{1}{4} \int_\Sigma |\vec{H}|_h^2$. There is an $\epsilon_0>0$ depending on $C_B$ and $K$ such that at any concentration point $p$ of $\HH$ around which the concentrating surfaces obey $\HH[\Sigma_r] \leq 4 \pi + \epsilon_0^2$,  we have
\[ \nabla^M \left(  \Sc_p + \frac{3}{5} \tr K_p ^2 + \frac{1}{5} |K_p|^2 \right) = 0 . \]
\end{theorem}

\begin{theorem}
\label{thm_intro_expansion}
Let $\Sigma \subset M$ be a spherical surfaces.  
Suppose $\Sigma$ is contained in a normal coordinate neighborhood $B_r(p)$ as in Lemma \ref{lemma_good_normal_coord} and that  $\| \A \|^2_{L^2(\Sigma)} \leq C r |\Sigma| $. Then $\mathcal{E}$ has the following expansion.
\[ \left| \mathcal{E}[\Sigma] - \frac{1}{12}\left(\frac{|\Sigma|}{4\pi}\right)^{3/2} \left(\Sc_p  + \frac{3}{5} \tr K_p^2 + \frac{1}{5} |K_p|^2 \right)  \right| \leq C |\Sigma|^2 \]
\end{theorem}
Note that Theorem \ref{thm_intro_expansion} stands in contrast to the results of G. Horowitz and B. Schmidt \cite{Horowitz_Schmidt_Note_on}. There they found that the Hawking energy has the following expansion 
\[\mathcal{E}[S] \sim \rho_p R^3  + O(R^4) \sim \left( \Sc_p + \tr K^2_p - |K_p|^2 \right) R^3 + O(R^4) ,\]
when calculated on spherical cross sections of the light cone in the tangent space at $p$.

This discrepancy is very surprising.
In general relativity the energy density, time component of the stress-energy tensor, $\rho$ is given by $ 16 \pi \rho = \Sc + (\tr K)^2 - |K|^2 $. 
 
As the Hawking energy should serve as a quasi local energy one might think that surfaces with maximal area constrained Hawking energy would tend to concentrate around critical points of the energy density $\rho$ which is not the case. 
Similarly, one would expect to find the energy density in the expansion of the Hawking energy. 

The fact that the expansion in a space like slice does not capture the energy density, where as the expansion along a light cone does, is especially vexing as the spheres in the light cone can be though of as lying in a space like slice, belonging to a different time, themselves.


\section{Minimizers of Hawking Type Functionals}
\label{section_Hawking_intro}

For the rest of the paper we work in a three dimensional Riemannian manifold $(M,g)$. The main objects of our study are immersed surfaces $\Sigma \subset M$ and we regard the immersion $\phi: S \to \Sigma$ as a parametrization. Additionally, we regard all involved functionals as being functionals on the immersions or the surfaces, interchangeably.

In \cite{Me_minimizers_gen} we investigated generalized Willmore functionals in detail.
In particular, we prove a compactness result for stratified surfaces which enables direct minimization in the class of bubble trees. Moreover, we show that critical points of generalized Willmore functionals are always smooth away from finitely many points.
Below we present the relevant definitions.

\begin{definition}[see {\cite[Definition 2.1]{Me_minimizers_gen}}, cf. {\cite[Definition 1 and 2]{Chen_Li_14}}]

Let $(S,\eta)$ be a Riemann surface and let $(M^n,g)$ be an n-dimensional, orientated Riemannian manifold which we assume to be isometrically embedded in some $\R^N$. 

\begin{enumerate}
\item For $k\in  \Z$ and $p \in [1,\infty]$ we define the Sobolev spaces as follows:
\[W^{k,p}(S,M) := \left\{ \phi\in  W^{k,p}(S, \R^N) \suchthat \phi(S) \subset M \textup{ a.e.} \right\}. \] 

\item
An element $\phi \in W^{2,2}(S,M)$ is called \emph{conformal immersion}, if $\phi$ is an immersion almost everywhere and if there is a function $e^{2\lambda}: S \to \R$, called the \emph{conformal factor} of $\phi$ such that
\[ \phi^* g = e^{2\lambda} \eta . \]

	\item  
	We say $\phi: S \to M$ is a \emph{branched conformal immersion} with finitely many branch points $ B \subset S$, if $\phi \in W^{2,2}_{loc}(S \setminus B ,M)$ is a conformal immersion and if for all $p \in B$ there is an open neighborhood $U_p$ and a constant $C$ such that
	\begin{align*}
	\int_{U_p \setminus \{p\}} 1+ |\vec{A}|^2 \, \d \mu_g  \leq C .
	\end{align*}
	
	\item Set 
	\begin{align*}
	\FF(S,M) &:= \{ \phi \in W^{2,2} (S,M) \mid \phi \textup{ is branched, conformal, immersion}\\
	& \hspace{0.8cm} \textup{with  branch points } B; \phi \in  W^{1,\infty}_{loc}(S\setminus B, M) \}
	\end{align*}
	and for $a>0$ define $\FF_a(S,M):= \{ \phi \in \FF(S,M) \mid |\phi(S)|= a \}$

\end{enumerate}
\end{definition}

 For an immersion $\phi \in \FF(S,M) $, $\Sigma:= \phi(S)$, we use $\WW[\Sigma]$ and $\WW[\phi]$ interchangeably. 
Moreover, at times we write $\AAA[\Sigma]$ or $\AAA[\phi]$ for the area $|\Sigma|$ in order to emphasize its role as a functional.

In \cite[Theorem 3.1]{Kuwert_Li_conf_immersions} E. Kuwert and Y. Li showed that branched conformal immersions can be extended to $W^{2,2}$ maps.

\begin{definition}[see {\cite[Definition 2.1]{Me_minimizers_gen}},  cf. {\cite[Definition 3]{Chen_Li_14}} ]
A compact connected metric space $(S,d)$ is called a \emph{stratified surface with singular points $P$}, if $P \subset S$ is a finite set such that:
\begin{enumerate}
	\item the regular part, $S \setminus P$, is a smooth Riemann surface without boundary. It carries a smooth metric $\eta$, whose induced distance function agrees with $d$.
 	\item  Moreover, for each $p\in P$ there is a $\delta > 0$ such that $B_\delta(p) \cap P = \{p\}$ and $B_\delta (p) \setminus \{p\} = \bigcup_{i=1}^{m(p)} \Omega_i$. Here $1<m(p) < \infty$ and the $\Omega_i$ are topological discs with one point removed. Additionally, we assume that $\eta$ can be extended to a smooth metric on each $\Omega_i \cup \{ p \}$.
\end{enumerate}
\end{definition}

By abuse of notation we usually denote a stratified surfaces as $S= \bigcup_i S^i$ and refer to Riemannian metrics on $S$ instead of on every $\overline{S}^i$.   

\begin{definition}[see {\cite[Definition 2.4]{Me_minimizers_gen}}] \
 \begin{enumerate}
\item 
Associate to every stratified surface $S = \bigcup_i S^i$ its \emph{dual graph}, where the vertices correspond to the components $S^i$ and two vertices are joined by an edge whenever the corresponding $S^i$  are joined by a singular point. 

	\item 
A stratified surface whose regular part consists of punctured spheres and whose dual graph is a simple tree is called a \emph{bubble tree}. 
The constituting spheres are called \emph{bubbles}.
\end{enumerate}
\end{definition}

\begin{definition}[see {\cite[Definition 2.5 and Definition 4.4]{Me_minimizers_gen}}] 

\ \begin{enumerate}
	\item 
Let $S$ be a stratified surface with $S \setminus P = \bigcup_{i=1}^m S^i $ and let $M$ be a manifold of dimension three or higher. 
For $k\in \N$ and $p \in [1,\infty]$ denote by $W^{k,p}(S, M)$ the continuous maps $\phi: S \to M$ for which all $\phi|_{S^i}$ extend to  maps in $W^{k,p}(\overline{S^i},M)$.\\
 Additionally, we say that $\phi: S \to M$ is a (branched) immersion if all extensions $\phi|_{\overline{S^i}}$ are (branched) immersions.

\item  Let $S = \bigcup_{i=1}^m S_i$ be a stratified surface and let $\phi:S \to M$ be a continuous map into a manifold $M$. We say $\phi$ is a \emph{haunted} immersion, if it is constant on some, but not all, components of $S$ and an immersion on the rest. A component $S_i$ is called a \emph{ghost} if $\phi|_{S_i}$ is constant, otherwise it is called \emph{regular}.

\item 
Any functional defined for immersed surfaces we extend componentwise to haunted, immersed stratified surface.
\end{enumerate}
\end{definition}

For simplicity we restate the definition of generalized Willmore functionals for stratified surfaces in a three dimensional ambient manifold.
\begin{definition}[see {\cite[Definition 2.7]{Me_minimizers_gen}}]
\label{def_gen_will}
Let $S$ be a stratified surface and let $(M,g)$  be an oriented three dimensional Manifold. 
\ \begin{enumerate}
	\item A branched conformal immersion $\phi \in \FF(S, M)$ is said to solve a \textit{generalized Willmore equation} (away from the branch points) if
	\begin{align}
	\label{eq_gen_will}
	\Delta H + H|\A|^2 + F(\phi, \nabla \phi, \nabla^2 \phi) = 0 ,
	\end{align}
	where $F$ is such that locally in conformal coordinates with $\lambda \in L_{loc}^\infty$ we have
	\begin{align*}
	e^{2\lambda}F(\phi, \nabla \phi, \nabla^2 \phi) &\in L^1 + H^{-1} \ \ \text{if } \phi \in \FF(S,M)  \\
	e^{2\lambda} F(\phi, \nabla \phi, \nabla^2 \phi) \nu &\in W^{k-1,l}, \ l= \frac{2p}{2+p} + \epsilon \ \ \text{if } \phi \in W^{k+2,p} , \ p>2 , \ k\geq 0
	\end{align*}
	for some $\epsilon >0$.
	
	\item A functional $ \HH$ on $\FF(S,M)$ is called an $a$-\textit{generalized Willmore functional}	if
\begin{enumerate}
	\item for any $\phi\in \FF_a(S,M)$ a bound $\HH[\phi] \leq \Lambda$ implies a bound on the Willmore energy $\WW[\phi] \leq C(\Lambda, a, M, \HH)$. 
	
	\item $\HH$ is bounded from below on $\FF_a(S,M)$.
	\item $\HH$ is invariant under diffeomorphisms of $S$.
	
	\item Let $\{\phi_k\}$ be a sequence in $\FF(S,M)$ with conformal factors $e^{2 \lambda_k}$. For any finite set $\mathfrak{S}\subset S$ the weak convergence $\phi_k \rightharpoonup  \phi$  in $W^{2,2}_{loc}(S\setminus \mathfrak{S},M)$ together with $ \|\lambda_k\|_{L^\infty(K)} \leq C_K$ for any $K \subset \subset S\setminus \mathfrak{S}$ implies $\HH[\phi] \leq \lim_{k\to \infty} \HH[\phi_k]$.	
	
	\item $\HH$ is differentiable and its Euler-Lagrange equation is a generalized Willmore equation.
\end{enumerate}
If a functional is an $a$-generalized Willmore functional for all $a>0$ or if the area in question is understood we will simply refer to it as a generalized Willmore functional.
\end{enumerate}
\end{definition}

As the Hawking energy was a prime motivator for this definition it will come to no surprise that Hawking type functionals are generalized Willmore.

\begin{prop}
Let $\HH_L$ be of Hawking type. Suppose that $L$ is smooth and bounded, then $\HH$ is a generalized Willmore functional.
\end{prop}

\begin{proof}
Let $\phi: S \to \Sigma \subset M$ be a closed, branched, immersed Riemann surface with area $a$ and $\HH[\Sigma] \leq \Lambda$. Suppose $|L|\leq C(L)$ then we have
\[ \HH[\Sigma] = \WW[\Sigma] + \int_\Sigma L(x,\nu) \,\d\mu > -C(L) a, \]
and
\[ \WW[\Sigma] \leq \Lambda - \int_\Sigma L(x,\nu) \,\d\mu \leq \Lambda + C(L) a.  \]

\noindent $\HH$ is invariant under reparametrisations of $\Sigma$ as $L$ is defined on $TM$. 

It is known that the Willmore energy is lower semi continuous in this setting (see for instance \cite[Lemma A.8]{Mon_Riv_13_1}). Thus we only need to discuss the lower order terms.

The convergence in $W^{2,2}_{loc}$, $\phi_k \to \phi$ implies local convergence in $W^{1,q}$ for all $1\leq q < \infty$. Thus we have point wise convergence almost everywhere of $\nabla \phi_k \to \nabla \phi$, $\phi_k \to \phi$ and hence of $\nu_k \to \nu$. Since $L$ is smooth, dominated convergence yields that $\HH$ is lower semi continuous. We examine the Euler-Lagrange equation of $\HH$ in the subsequent lemmas.
\end{proof}

Let $\phi_0 : S \to \Sigma \subset M$ be a $W^{2,2}(S,M)$ immersed Riemann surface. 
Consider a normal variation of $\phi_0$ along the vector field $f \nu$, $f \in C^\infty(\Sigma)$, that is
\begin{align*}
\phi: (-\epsilon, \epsilon)\times S & \to M \\      
(s,x) &\mapsto \phi(s,x) = \phi_s(x)
\end{align*}
 such that for every $s\in (-\epsilon, \epsilon)$, $\phi_s(S)$ is an immersed surface in $M$, $\phi(0,S) = \phi_0(S)= \Sigma$ and $\left. \frac{\partial \phi}{\partial s}\right|_{s=0} = f \nu$.

\begin{lemma}
\label{lemma_basic_variation}
Under the variation above the geometric quantities behave as follows.   
\begin{align}
\frac{\partial}{\partial s}\Big|_0 d \mu & =  f H d\mu \nonumber \\ 
\frac{\partial}{\partial s}\Big|_0\nu & =  - \nabla^\Sigma f  \nonumber \\              
\frac{\partial }{\partial s}\Big|_0 H & =  - \Delta f - f\left(  |A|^2 + \Ric^M (\nu,\nu) \right) \nonumber \\ 
\frac{\partial  }{\partial s}\Big|_0 L(\phi_s, \nu_s)  &= d_{TM}L ( f \nu, - \nabla^\Sigma f) \nonumber \\
\label{eq_var_L}
&= f d_M L (\nu) - d_V L (\nabla^\Sigma f)
\end{align}
Here $d_{TM}$ denotes the exterior differential of $TM$, and  $d_M h(X) := d_{TM}h (X,0)$ and $d_V h(X) :=d_{TM}h(0,X) $ for a vector field $X$ on $M$ and  $h \in C^1(TM)$.
The corresponding area constrained Euler-Lagrange equation  reads
\begin{align} 
\label{eq_EL_H}
\Delta H + H |\A|^2 + H Q + \gamma (\A , S ) + 2 \lambda H + T =0 .
\end{align}
If $L$ is smooth then it is a generalized Willmore equation. Here $Q$, $S$ and $T$ are defined as
\begin{align*}
Q &:= \Ric^M(\nu,\nu) -2 L  -  \tr_\Sigma \Hess_V L + 2 d_V L (\nu)\\
S &:= - 2 \Hess_V L \\
T &:= - 2 d_M L(\nu) - 2 \Div_\Sigma d_V L,
\end{align*}
where $\Hess_V L$ is the fiber part of the Hessian of $L$ and $\Div_\Sigma d_V L = \tr_\Sigma \nabla^M d_VL$ with\\ $\nabla^M d_V L (X,Y) = \nabla^{TM}_{(X,0)} d_{TM}L (0,Y)$.
\end{lemma}       

\begin{proof}
The variation of the geometric quantities is widely known, see  for instance \cite[Theorem 3.2, Section 7]{Huisken_Polden_Geo_Evo}  and the variation of $L$ is straight forward.

The variational problem reads $\lambda \delta_f \AAA(\Sigma) = \delta_f \WW(\Sigma) + \delta_f \LL(\Sigma) $.
We treat all terms separately.
\begin{align*}
\frac{d}{d s}\Big|_0 \mathcal{A}[\Sigma] &= \int_{\Sigma} f H  \, \d \mu  \\
 \frac{d}{d s}\Big|_0 \WW[\Sigma] &= \frac{1}{2} \int_\Sigma  H \frac{\partial H}{\partial s}\Big|_0 + f \frac{H^3}{2}  \, \d\mu \\
&=\frac{1}{2} \int_\Sigma  H  (- \Delta f - f\left(  |\A|^2 + \Ric^M (\nu,\nu) \right))  \,  \d\mu\\
&= \frac{1}{2} \int_\Sigma - f \left( \Delta H +  H|\A|^2 + H \Ric^M (\nu,\nu) \right)  \, \d\mu \\
 \frac{d}{d s}\Big|_0 \LL[\Sigma] &= \int_\Sigma \frac{\partial L}{\partial s}\Big|_0 + f L H \, \d \mu \\
\int_\Sigma d_V L (\nabla^\Sigma f) \, \d\mu &= \int_\Sigma   \Div_\Sigma (f d_V L) - f \Div_\Sigma \left( d_V L_{(x,\nu)} \right)  \, \d \mu \\
&= - \int_\Sigma f \left( (\Div_\Sigma d_V L)_{(x,\nu)} + \gamma ((\Hess_V L )_{(x,\nu)}, A) -H d_V L(\nu)   \right) \, \d \mu 
\end{align*} 
Sorting all the terms yields the desired equation.
\[ \Delta H + H |\A|^2 + H Q + \gamma (\A , S ) + 2 \lambda H + T =0  \]
In the notation of generalized Willmore equations we have $F =  H Q + \gamma (\A , S ) + 2 \lambda H + T  \in L^{2} $, provided $\phi \in \FF(S,M)$.

For the higher order regularity note that the worst term of $|\nabla \phi|^2 F \nu$ is of the form $|\nabla \phi|^2 \star A \star \psi(\phi,\nu)  \star \nu\star \nu \star \nu$, where the $ \star $ denotes a sum of contractions, and $\psi:TM\to\R$ is a smooth and bounded function. If $\phi \in W^{k+2,p}\cap W^{1,\infty}$, $k\geq 0$, $p\geq 2$, then, due to the Sobolev embeddings $W^{k+2,p} \hookrightarrow W^{k+1,q} \hookrightarrow C^{k,\alpha}$, for all $1\leq q< \infty$ and $\alpha \in (0,1)$, we have $|\nabla \phi|^2  \star \nu\star \nu \star \nu \in W^{k+1,p} \cap L^\infty$. Similarly, $\psi(\phi,\nu) \in W^{k+1,p} \cap  L^\infty$. This means, due to the $A$ component we have $|\nabla \phi|^2 F \nu \in W^{k,p}$, whenever $\phi \in W^{k+2,p}\cap W^{1,\infty}$.
\end{proof}

\begin{lemma}
\label{lemma_con_check_1}
Let $X$ and $Y$ be a vector fields along $\Sigma$, and introduce the 1-form $\eta := K(\cdot,\nu)$ as well as the musical isomorphism $\#: T^*\Sigma \to T\Sigma$  then the following equations hold.
\begin{align*}
d_M P^2 (X) &= 2 P \tr_\Sigma \nabla^M_X K \\
d_V P^2 (X) &= -4 P K(X, \nu) \\
\nabla^M_Y d_V P^2 (X) &= -4 \left( tr_\Sigma \nabla^M_Y K \right) K(X,\nu) - 4 P \nabla^M_YK(X,\nu)\\
\Hess_V P^2 &= 8 \eta \otimes \eta - 4 P K \\
\Div_\Sigma d_V P^2 &= - 4 \tr_\Sigma \nabla^M_{\eta^\#} K - 4 P \Div_\Sigma K (\nu) 
\end{align*}

Moreover, the area constrained Euler-Lagrange equation for the Hawking type functional with $L= - \frac{1}{4}P^2$ reads
\[ \Delta H + H |\A|^2 + H Q + f \gamma( \A , S ) + 2 \lambda H + T =0 . \]
  $Q$, $S$ and $T$ are given by
\begin{align*}
Q &= \Ric^M(\nu,\nu) -\frac{1}{2} P^2  + 2|\eta|^2 + 2P K (\nu,\nu)\\
S &= -2P K + 4 \eta\otimes \eta\\
T &=  P \tr_\Sigma \nabla_\nu K  - 2 \tr_\Sigma \nabla_{\eta^\#} K - 2 P \Div_\Sigma K(\nu) .
\end{align*}
\end{lemma}

\begin{proof}
The first equation is clear because the metric is parallel. 
For the second we use $P(x,\nu_x)= \tr K_x - K_x(\nu_x,\nu_x) = \left( \tr K - K \right)_{(x,\nu_x)} $ and find that 
\[ d_V P (X) = - 2K(X,\nu) \]
as $K$ is bilinear and $\tr K$ does not depend on the fiber.
Now the rest follows easily.
\begin{align*}
d_V P^2 (X) &= -4 P K(X,\nu)\\
\nabla^M_Y d_V P^2 (X) &= -4 \left( \tr_\Sigma \nabla^M_Y K \right) K(X,\nu) - 4 P \nabla^M_YK(X,\nu)\\
\Hess_V P^2 &= 2 d_V ( -2 P K(\cdot,  \nu) ) \\
&= 8 K(\cdot, \nu) \otimes K(\cdot, \nu) -4 P K \\
\Div_\Sigma d_V P^2 &= -4 \Div_\Sigma \left( P K(\cdot, \nu) \right) \\
&= -4 P \left( \Div_\Sigma K \right)(\nu) - 4 \gamma^{ij} \tr_\Sigma (\nabla^M_i K) \eta_j \\
&= -4 P \left( \Div_\Sigma K \right)(\nu) - 4\tr_\Sigma \nabla^M_{\eta^\#} K
\end{align*}
The Euler-Lagrange equation is obtained from \eqref{eq_EL_H}.
\end{proof}

Next, we will construct minimizers of general Willmore functionals through direct minimization.
Denote by $\TT $ the class of bubble trees, and define 
\[\FF_a(\TT,M) := \left\{ \phi \in \FF(S,M) \mid S \in \TT, \ (S,\phi) \text{ haunted, } \AAA[\phi] =a \right\} \]
as well as 
\[ \beta(\HH, M, a) = \inf \{\HH[\phi] \mid \phi \in \FF_a(\TT,M) \} . \]

\begin{theorem}
\label{thm_small_minimizers}
Let $(M,g)$ be $C_B$-bounded (see Appendix \ref{chapter_bounded_geo_small_surf}) and let  $\HH_L$ be a Hawking type functional for a smooth and bounded $L$.
\begin{enumerate}
	\item If $M$ is compact, then there exits a $\phi \in \FF_a(\TT, M) $ realizing $\beta(\HH_L, M, a)$, for any $a >0$. Additionally, $\phi:S \to M $ is smooth away from its finitely many branch points and the finitely many singular points of $S$.
\item
There are constants $a_0(L, C_B)>0$ and $C(L, C_B)>0$ such that any $\phi_a \in \FF_a(\TT, M) $  realizing  $\beta(\HH_L, M, a)$ for $a\leq a_0$ is an embedding of a sphere, its image is contained in a normal coordinate neighborhood and $\phi_a$ satisfies
 \[ | \HH[\phi_a] - 4\pi| \leq C(L, C_B) a.\]
\end{enumerate}
\end{theorem}

\begin{proof}
Since we know that $H_L$ is a generalized Willmore functional, the first statement follows directly from \cite[Theorem 4.5 and Theorem 5.6 ]{Me_minimizers_gen}. 

 In \cite{Mondino2010} A. Mondino calculated the expansion of the Willmore energy for spheres in coordinates and found 
\[\WW[S_R, g]  \leq 4\pi + C(C_B) |S_R|_g. \]
 Since $L$ is bounded, we can estimate the Willmore energy of $\Sigma_a$ by comparing it to spheres $S_R$, $a= 4\pi R^2$, in coordinates.
\[ \WW[\Sigma_a,g] \leq \HH[\Sigma_a] + C(L) a \leq \HH[S_R] + C(L) a \leq 4 \pi + C(L, C_B) a \] 
Lemma \ref{lemma_diameter} asserts that $\diam_M (\Sigma_a) \leq C(L, C_B) \sqrt{a}$. 
Hence, $\Sigma_a$ lies in some normal coordinate neighborhood $B_r(p)$, provided $a$ is small enough. \\
Consider a bubble $S^i$ on which $\phi$ is not constant and set $\Sigma^i := \phi(S^i)$.
We apply corollary \ref{cor_small_surf_euc} as well as the integrated Gauss equation to see 
\[ \WW[\Sigma^i, g] \geq \WW[\Sigma^i, g_E] - C r^2 \geq 4 \pi - C r^2 . \]
Hence, there can be only one bubble and by deleting ghosts we may assume that $(S,\phi)$ is not haunted. From the same corollary we gather
\[ \WW[\Sigma_a, g_E] \leq \WW[\Sigma_a,g] + C r^2 \leq 4 \pi + C(r^2 +a ). \]
In order to see embeddedness we employ the Li-Yau inequality \cite{Li-Yau-inequality}. Denote by $\theta^2(\Sigma, p) = \# \phi^{-1}(p)$ the density of $\Sigma$ at $p$, then we have
\[ \theta^2(\Sigma, p) \leq \frac{\WW[\Sigma]}{4\pi}. \]
This follows from Simons monotonicity formula, see \cite[Appendix A]{Kuwert_Schaetzle_removability} for a discussion.

Finally, Lemma \ref{lemma_good_normal_coord} allows us to choose the normal neighborhood such that $r$ and $R$ are comparable. This yields the final estimate on $\HH[\Sigma_a]$.
\begin{align*}  
\HH[\Sigma_a]  & =  \WW[\Sigma_a, g] + \LL[\Sigma_a] \leq 4\pi + C(L,C_B) a \\
\HH[\Sigma_a] & \geq \WW[\Sigma_a, g_E] - C(L, C_B) a \geq 4 \pi - C(L, C_B) a
\end{align*}

\end{proof}


\section{Concentration of Critical Surfaces}
\label{section_concentration}

In this section we analyze critical points of a generalized Willmore functional of Hawking type with small area. We follow \cite{LMII} closely, where these arguments were developed for the Willmore functional. To that end we fix a $C_B$-bounded three dimensional ambient manifold $(M,g)$ (see Appendix \ref{chapter_bounded_geo_small_surf}) and a Hawking type functional $\HH = \WW + \LL$, where $\LL[\Sigma] = \int_\Sigma L(x,\nu) \, \d \mu$, for a smooth $L:TM \to \R$.
Throughout Subsection \ref{section_a_priori_small_surf} we suppose also that $L$, $d_{TM}L$, $\Hess_{V}L$ and $\nabla^M d_V L$ are bounded by $C_L$.
 As a general naming scheme we adopt that quantities on a surface $\Sigma$ are denoted by their usual symbol, whereas quantities on $M$ get identifying indexes. Geodesic balls in $M$  will be denoted by $\BB_r(p)$.

\subsection{A Priori Estimates for Small Critical Surfaces}
\label{section_a_priori_small_surf}

\begin{prop}
There are positive constants $r_0(C_B)$ and $C(C_B, C_L)$ such that for all $r \in (0,r_0)$ and $\Sigma \subset \BB_{r}(p)$, immersed, area constrained, critical surfaces of $\HH$, we can estimate the Lagrange multiplier as follows.
\[ |\lambda| \leq C |\Sigma|^{-1} \left(|\Sigma|^{1/2}\WW[\Sigma]^{1/2} + |\Sigma| +  r \int_\Sigma |A|^2 \d\mu \right) \]
\end{prop}

\begin{proof}
As in  \cite[Proposition 5.3]{LMII} the idea is to consider an area constrained normal variation of $\HH$ in direction $f = g(x,\nu)$, for the position vector field $x$ in $B_r$ in normal coordinates to obtain
\[ \delta_f \HH[\Sigma] = \lambda \delta_f \AAA[\Sigma]. \]
If the variation of the area is non zero, we calculate the Lagrange parameter $\lambda$ as the quotient
\[\lambda =   \frac{\delta_f \HH[\Sigma]}{\delta_f \AAA[\Sigma]}. \]
 From \cite{LMII} we know
\begin{align*}
\delta_f \AAA[\Sigma] &\geq |\Sigma| \\
|\delta_f \WW[\Sigma]| & \leq C(C_B) |\Sigma|^{1/2} \WW[\Sigma]^{1/2} + C(C_B) r(1+|\Sigma|^{1/2}) \int_\Sigma |A|^2 \,\d\mu,
\end{align*}
and we estimate $\delta_f \LL$ as follows.
\begin{align*}
| \nabla^\Sigma f | &\leq \sum_{i} | \gamma^{ij} g(\nabla^M_j x , \nu ) | + \left| \gamma^{ij}  g\left(x, \frac{\partial \phi}{\partial x^k}\right) A^k_j  \right| \\
&\leq C(C_B) + C(C_B) r |A| \\ 
 \delta_f \LL[\Sigma]  &=  \int_\Sigma f d_M L (\nu)+ d_V L(-\nabla^\Sigma f) + fL H  \, \d \mu \\
& \leq C(C_L,C_B)|\Sigma| + C(C_L,C_B) r |\Sigma|^{1/2} \left( \int_\Sigma |A|^2 \,\d\mu \right)^{1/2} 
\end{align*}
\end{proof}

\begin{prop}
\label{thm_a_priori}
There are positive constants $\epsilon_0 $ and $C$ depending only on $C_B$ and $C_L$ such that any spherical immersed surface $\Sigma$ that 
\begin{enumerate}
\item solves equation \eqref{eq_EL_H}, satisfies
\item $\HH(\Sigma) \leq 4 \pi + \epsilon^2$ and
\item $|\Sigma|\leq \epsilon^2 $
\end{enumerate}
for an $\epsilon \in (0,\epsilon_0)$, obeys the following estimate.
\[ \int_\Sigma | \nabla^2 H |^2 + H^2 |\nabla H|^2 + H^2|\nabla \A|^2 + H^4 |\A|^2 d \mu \leq C\] 
\end{prop}

\begin{proof}
We start by integrating the Gauss equation over $\Sigma$, to obtain
\[ 2 \WW[\Sigma] = 8 \pi + \int_\Sigma |\A|^2 d\mu + 2 \int_\Sigma G^M(\nu,\nu)  \, \d \mu .\]
Here $G^M = \Ric^M - \frac{1}{2} Sc^M g$ denotes the Einstein tensor of $M$. Since $G^M$ is bounded and $\HH$ is close to $4 \pi$, we can estimate $\| \A \|^2_{L^2(\Sigma)}$.
\begin{align*} \| \A \|^2_{L^2(\Sigma)} &= 2 \HH[\Sigma] - 8 \pi -2 \LL[\Sigma] -2 \int_\Sigma G^M(\nu,\nu)\, \d \mu \\
& \leq C(C_L, C_B) \epsilon^2
\end{align*}
Moreover, Lemma \ref{lemma_diameter} and Lemma \ref{lemma_good_normal_coord} assert that we can operate in a normal coordinate neighborhood $\BB_r(p)$ adapted to $\Sigma$ such that $r \leq C |\Sigma|^{1/2}$. 
This simplifies the estimate for the Lagrange multiplier to 
\begin{align*}
|\lambda|\leq C(C_L,C_B) |\Sigma|^{-1} \epsilon .
\end{align*}
We multiply equation \eqref{eq_EL_H} by $\Delta H$ and integrate over $\Sigma$. Through integration by parts and Young's inequality we obtain 
\begin{align*}
\int_\Sigma (\Delta H)^2 \, d\mu &= - \int_\Sigma \Delta H H |\A|^2 + \Delta H H Q + \Delta H \gamma(\A, S) + \Delta H T -2 |\nabla H|^2 \lambda \, \d \mu  \\
& \leq \int_\Sigma \frac{(\Delta H)^2 }{2} + 2 H^2 |\A|^4 + 2 H^2 Q^2  + 2  \gamma(\A, S)^2 + 2 T^2 +\frac{C \epsilon}{|\Sigma|} |\nabla H|^2 \, \d \mu  ,
\end{align*}
and
\begin{align} 
\int_\Sigma (\Delta H)^2  \, d\mu & \leq C(C_L, C_B)  + C(C_L,C_B) \frac{\epsilon}{ |\Sigma|} \int_\Sigma |\nabla H|^2 \, d\mu  + 4 \int_\Sigma  H^2 |\A|^4 \, \d \mu  .
\label{eq1_part1}
\end{align}
From here on the integral estimates are identical to the proof of \cite[Proposition 5.1]{LMII}.

\end{proof}

The next corollary establishes the roundness of small surfaces of generalized Willmore type. It is virtually the same as \cite[Corollary 5.5]{LMII}. We will not prove it here as it relies only on the estimate of Proposition \ref{thm_a_priori}, which is the analog of \cite[Theorem 5.4]{LMII}, and general facts about small surfaces.

\begin{corollary}
\label{cor_H_L_infty}
Assume $\Sigma$ is a surface as in Proposition \ref{thm_a_priori}.  If $|\Sigma|$ is small enough, there exists a constant $C = C(C_L, C_B)$ such that the following estimates hold.
\begin{align*}
 \|\A\|_{L^2(\Sigma)} & \leq C |\Sigma| \\
\|H- 2/R \|_{L^\infty(\Sigma)} &\leq C |\Sigma|^{1/2} 
\end{align*}
In particular, the mean curvature has to be positive and the inverse of the mean curvature has to be bounded. 
\begin{align*}
\| H^{-1} \|_{L^\infty(\Sigma)}  &\leq C |\Sigma|^{1/2}\\
\| H^{-1} - R/2  \|_{L^\infty(\Sigma)} &\leq C |\Sigma|^{3/2}
\end{align*}
\end{corollary}


\subsection{Surface Concentration}
\label{section_surface_concentration}

In this section we characterize the points around which surfaces of generalized Willmore type concentrate. 
It is a direct generalization of the corresponding results obtained by T.Lamm and J.Metzger in \cite{LMI}.

\begin{definition}
A Point $p\in M$ is called a concentration point of $\HH$ if there is a constants $r_0 > 0$ and an $A_0 > 0$ such that for every $r\in (0,r_0)$ there is an $A \in (0,A_0)$ and a spherical, area constrained, critical surface $\Sigma_r$ of $\HH$ with $|\Sigma_r| = A$ contained in the geodesic ball $\BB_r(p)$. 
\end{definition}

\begin{definition}
Let $S:= S^2_1(a)$ be the two sphere around $a\in \R^3$ with outer normal vector field $\nu$ and fet $F: \R^3 \times \R^3 \to \R$ be bounded. 
For a multi index $(\alpha_1,...,\alpha_k)$, $k\in \N$ introduce 
\begin{align*}
c^{(\alpha_1,...,\alpha_k)}(F, a) &:= \int_{S^2_1(a)} F(a,\nu) (x-a)^{\alpha_1} ...(x-a)^{\alpha_k} \, \d\m , \\
c(F, a) &:= \int_{S^2_1(a)} F(a,\nu) \, \d\mu .
\end{align*}
\end{definition}

\begin{theorem} 
\label{thm_concentration_general}
Let $(M,g)$ be a $C_B$ bounded three manifold and let $\HH = \WW + \LL$ be a Hawking type functional with $\LL[\Sigma] =  \int_\Sigma  L(x,\nu) \, \d \mu$ for a smooth $L$.

\begin{enumerate} 
  \item Let $M$ be compact, then there exits at least one concentration point of $\HH$. The concentrating surfaces $\Sigma_r$ at that point are area constrained minimizers of $\HH$ and obey $\HH[\Sigma_r] \leq 4 \pi + \epsilon_0^2$, where $\epsilon_0$ is the constant from Proposition \ref{thm_a_priori}.

	\item 
Let $p \in M$ be a concentration point of $\HH$ such that the concentrating surfaces $\Sigma_r$ have energy $\HH[\Sigma_r] \leq 4 \pi +\epsilon_0^2$.\\
Then in Riemannian normal coordinates around $p$ the vector $V_p$ with components 
\[ V^\alpha_p = -c \left( d_VL_\alpha,p \right) + 2 c^\alpha \left( L,p  \right)  +  c^{(\alpha,\beta)} \left( d_VL_\beta,p \right) \] 
vanishes.\\
Moreover, if $V_{a_i}$ vanishes identically for a sequence of points $\{a_i\}$ converging to $p$, as constructed in the proof, then we have that
\[  \nabla^M \Sc_p  - W_p  = 0 .  \]
Here $W_p$ is a vector whose components read 
\[W^\alpha_p =  \frac{3}{2 \pi} \left(  - c^\beta \left(\nabla^M_\beta d_VL_\alpha, p \right)  + 3 c^{(\alpha,\beta)}\left( d_ML_\beta,p \right) + c^{(\alpha,\beta,\gamma)}\left( \nabla^M_\gamma d_V L_\beta ,p \right)  \right) . \]
\end{enumerate}

\end{theorem}

\begin{remark} \ 
\begin{enumerate}
\item[a)]
If $L(x,\nu)$ is even in $\nu$, then the $V_a$ vanish as all the involved integrals vanish.

\item[b)] $W_p$ involves only terms with a $d_M$. We can therefore see it as the gradient of some function $w$ at $p$. This leads to the interpretation that, provided $V$ vanishes, a concentration point of $\HH$ is a critical point for $\Sc - w$.
\end{enumerate} 
\end{remark}

\begin{proof}
For the first part, we know by Theorem \ref{thm_small_minimizers} that there is a minimizing area constrained embedded sphere with $\HH[\Sigma_A] \leq 4 \pi + \epsilon_0^2$, for any small enough area $A$. Moreover, they are contained in normal neighborhoods $\BB_{r_A}(p_A)$, where $r_A$ and $\sqrt{A}$ are comparable. 
For $A \to 0$ the points $p_A$ will subconverge to a point $p$ which is a concentration point by construction.

For the second part, let $r_0$ and $A_0$ be as in the definition of concentration point.
Suppose $r \in (0,r_0)$ and  $r_0 \leq \inj(M,g)$. 
Let $\Sigma$ be a spherical, area constrained, critical point of $\HH$ contained in  $\BB_{r}(p)$ with area $|\Sigma| = 4 \pi R^2$ and $\HH[\Sigma_r] \leq 4 \pi +\epsilon_0^2$. 
Since $L $ is smooth and we work in $\BB_r(p)$, the results of Section \ref{section_a_priori_small_surf}
apply. 
In Appendix \ref{chapter_bounded_geo_small_surf} we discuss that, by choosing $r_0$ smaller if necessary, we have the estimates $d:= \diam_g \Sigma \leq C R \leq C r$. 
Since there is at least one such $\Sigma$ for any $r\in (0,r_0)$, we may suppose that $2d < r$. 
This allows us to use Lemma \ref{lemma_good_normal_coord} to find normal coordinates $\psi$ adapted to $\Sigma$ around $p_\Sigma \in M$ such that $\Sigma \subset \BB_{2 d} (p_\Sigma) \subset \BB_{r}(p_\Sigma)$, $d_g(p,p_\Sigma) \leq r$ and
\[ \int_{\psi(\Sigma)} y \, \d \mu_g(y) = 0 . \]
Additionally, in these adapted normal coordinates we have  
\[\max_{x \in \Sigma}|x|_E \leq C R. \]
We will operate in these coordinates from now on.

Consider the area constrained variation of $\HH[\Sigma]$ with respect to the vector field $f \nu = g(b,\nu)/H \nu$, where $b$ is a constant vector field to be chosen later.
Recalling the traced Gauss equation
\[ \Sc^\Sigma = \Sc^M - 2 \Ric^M(\nu,\nu) + \frac{1}{2} H^2 - |\A|^2,\]
we may split the Willmore functional into two new functionals
\begin{align*}
\UU[\Sigma] &= \frac{1}{2}\int_\Sigma |\A|^2 \, d\mu \\
\VV[\Sigma] &= \int_\Sigma  \Ric^M(\nu,\nu) - \frac{1}{2}\Sc^M \,d\mu \\
\WW[\Sigma] &=4 \pi(1-q(\Sigma))+ \UU[\Sigma] +  \VV[\Sigma]
\end{align*}
and arrive at
\begin{align*}
\lambda \delta_f \AAA =  \delta_f \HH =  \delta_f \UU + \delta_f \VV + \delta_f \LL .
\end{align*}
Let $\Omega$ be the region enclosed by $\Sigma$ and let $\langle \cdot , \cdot \rangle$ denote the Euclidean scalar product on $\R^3$.  

Estimating $\Vol(\Omega)$, $\delta_f\AAA[\Sigma]$ and $\delta_f\UU[\Sigma]$, as well as the better part of $\delta_f\VV[\Sigma]$ as in \cite[Section 4]{LMI} yields
\begin{align*}
\left| \Vol(\Omega) - \frac{4 \pi}{3} R^3 \right| & \leq C R^5\\
|\lambda \delta_f \AAA[\Sigma]| &\leq C R^4 \\
|\delta_f \UU[\Sigma]| &\leq C R^4 \\
\left|\delta_f \VV[\Sigma] + \frac{1}{4} \Vol(\Omega) \langle \nabla^M \Sc_{p_\Sigma} , b \rangle \right| &\leq C R^4 \\
\end{align*}
and hence
\begin{align}
\label{eq_crit_expansion}
\left| \frac{\pi}{3} R^3 \langle b, \nabla^M \Sc^M_{p_\Sigma}  \rangle - \delta_f \LL[\Sigma] \right| \leq C R^4 .
\end{align}
Thus we need to estimate the variation of $\LL$.
\begin{align}
\label{eq_var_LL_1}
\delta_f \LL[\Sigma] &= \int_\Sigma f d_M L (\nu) - d_V L(\nabla f) + L H f \, d \mu
\end{align}
We start with the second term on the right hand side, using $e_i := \frac{\partial \phi}{ \partial x_i}$.
\begin{align*}
\int_\Sigma  d_V L (\nabla f)   \, d\mu &= \int_\Sigma \frac{1}{H} d_VL(e_j) \gamma^{ij} \left( g(\nabla^M_i b, \nu) + A(b^T, e_i) - \frac{1}{H} \frac{\partial H}{ \partial x^i} g(b,\nu) \right) \, d\mu \\
&= \int_\Sigma   \frac{1}{H} d_VL(e_j) \gamma^{ij}   \left( \A(b^T,e_i) + g(\nabla^M_i b, \nu) \right) - \frac{1}{H^2} d_VL( \nabla H)   \\
  & \hspace{1cm}+ \frac{1}{2} d_VL(b^T)  \,d\mu 
\end{align*}
The first three terms on the right hand side can be estimated rather easily, using the results of  Section \ref{section_a_priori_small_surf} as well as the fact that we use Riemann normal coordinates on $M$.
\[ \int_\Sigma   \frac{1}{H} d_VL(e_j) \gamma^{ij}   \left( \A(b^T,e_i) + g(\nabla^M_i b, \nu) \right) - \frac{1}{H^2} d_VL( \nabla H)  \,d\mu \leq C R^4 \]
The other terms, that is $\int_\Sigma 1/2 \, d_VL(b^T) + f d_ML(\nu) + L H f \,d\mu$, need to be treated in more detail. We will pull them back to an approximating sphere to perform explicit calculations. 
In Lemma \ref{lemma_small_surf_euc}, Theorem \ref{thm_sigma_to_s} and corollary \ref{cor_sigma_to_s} we detailed how this is possible. 
The estimates derived there in conjunction with corollary \ref{cor_H_L_infty} imply that, up to order $O(R^4)$ we have to estimate

\[ \int_S -\frac{1}{2} \, d_VL(b^T) +  \frac{R_E}{2} \langle b,\nu\rangle d_ML(\nu) +\langle b,\nu\rangle L  \,\d\mu .\]
Here  $R_E$ is the Euclidean radius of $\Sigma$, $|\Sigma|_E = 4\pi R_E^2$, which is comparable to $R$ and 
$S$ is the round sphere of radius $R_E$, centered at $a = \int_\Sigma x \, \d \mu_E$, the Euclidean center of mass of $\Sigma$. 
The outer normal to $S$ is given by $\nu = (x-a)/R_E$, where $x$ is the position vector field.\\
Note that in this construction $d_g(p_\Sigma,a) \leq C d^3$ and hence 
\begin{align}
\label{eq_p_a}
d_g(p, a) \leq C r .
\end{align}
Note  also that the term $1/2 d_VL(b^T) + \langle b,\nu \rangle L$ is of order one, whereas $R_E/2 \langle b,\nu \rangle d_M L(\nu)$ is of order $R$. This means, unless $\int_\Sigma 1/2 d_VL(b^T) + \langle b,\nu \rangle L$ vanishes up to $O(R^2)$, it will dominate the concentration point $p$. We will perform a Taylor expansion in the first variable around $a$ in order to separate the orders of magnitude. 
\begin{align*} 
d_VL(b^T)_{(x,\nu)} &= d_VL(b^T)_{(a,\nu)} + \nabla^M_{x-a} d_VL(b^T)_{(a,\nu)} + O(R^2) \\
\langle b,\nu \rangle L(x,\nu) &=\langle b,\nu \rangle L(a,\nu) + \langle b,\nu \rangle d_M L(x-a)_{(a,\nu)} + O(R^2) \\
\frac{R_E}{2} \langle b,\nu \rangle d_ML(\nu)_{(x,\nu)} &= \frac{R_E}{2} \langle b,\nu \rangle d_ML(\nu)_{(a,\nu)} + O(R^2)
\end{align*} 
Integrating over $S$ and separating by powers of $R_E$ yields
\begin{align*}
\int_S - \frac{1}{2} \,&  d_VL(b^T)  +  \frac{R_E}{2} \langle b,\nu \rangle d_ML(\nu) + \langle b,\nu \rangle  L  \,\d\mu \\
& = R_E^2\left(  -\frac{b^\alpha}{2} c(d_VL_\alpha) + b^\alpha c^\alpha(L)  + \frac{b^\alpha}{2} c^{(\alpha,\beta)}(d_VL_\beta) \right) \\
& \hspace{0.2cm} + R_E^3\left(  - \frac{b^\alpha}{2} c^\beta(\nabla^M_\beta d_VL_\alpha)  + \frac{3 b^\alpha}{2} c^{(\alpha,\beta)}(d_ML_\beta) + \frac{b^\alpha}{2} c^{(\alpha,\beta,\gamma)}(\nabla^M_\gamma d_VL_\beta) \right) \\
& \hspace{0.2cm} + O(R^4) .
\end{align*}
Define the components of two vectors $V_a$ and $W_a$ by 
\begin{align*}
V^\alpha_a &= (-c(d_VL_\alpha) + 2 c^\alpha(L)  +  c^{(\alpha,\beta)}(d_VL_\beta)) \\
W^\alpha_a &=  \frac{3}{2 \pi} \left(  - c^\beta(\nabla^M_\beta d_VL_\alpha)  + 3 c^{(\alpha,\beta)}(d_ML_\beta) + c^{(\alpha,\beta,\gamma)}(\nabla^M_\gamma d_VL_\beta )  \right) .
\end{align*} 
If $V_a$ is not  zero, then equation \eqref{eq_crit_expansion} implies that $\langle b, V_a \rangle  \to 0 $ for $r\to 0$ and any constant vector $b$. Moreover, by equation \eqref{eq_p_a} we get that $a \to p$ as $r \to 0$. Choosing $b= V_p$ yields that $p$ is characterized by the vanishing of $V_p$.\\
If $V_a$ vanishes, we get
\[  \nabla^M \Sc_p  - W_p  = 0,  \]
using equation \eqref{eq_crit_expansion} and $b = \nabla^M \Sc_p  - W_p $.

\end{proof}

Now we apply the previous result to the Hawking energy $\mathcal{E}$. Recall 
\[\mathcal{E}[\Sigma] = \frac{|\Sigma|^{1/2}}{16 \pi^{3/2}} \left( 4\pi - \HH[\Sigma] \right) \]  for $\HH[\Sigma] = \WW[\Sigma] - \frac{1}{4} \int_\Sigma (\tr_\Sigma K)^2 \, \d \mu $, where $K$ is a smooth symmetric 2-tensor on $M$. Clearly the area constrained minimizers of $\HH$ are the area constrained maximizers of $\mathcal{E}$.

\begin{theorem} 
\label{theorem_Hawking_concentration}
Let $\HH$ be as above. At any concentration point $p$ of $\HH$ around which the concentrating surfaces obey $\HH[\Sigma_r] \leq 4 \pi + \epsilon_0^2$, where $\epsilon_0$ is the constant from Proposition \ref{thm_a_priori}  we have
\[ \nabla^M \left(  \Sc_p + \frac{3}{5} \tr_M K_p ^2 + \frac{1}{5} |K_p|^2 \right) = 0 . \]
\end{theorem}

\begin{proof}
We apply Theorem \ref{thm_concentration_general}. First note that the function $L= -\frac{1}{4} (\tr_\Sigma K)^2$ is even in $\nu$ and hence the vectors $V_a$  vanish. Thus we need to compute the vector $W_p$ with components  
\[W^\alpha_p =  \frac{3}{2 \pi} \left(  - c^\beta(\nabla^M_\beta d_VL_\alpha,p)  + 3 c^{(\alpha,\beta)}(d_ML_\beta, p) + c^{(\alpha,\beta,\gamma)}(\nabla^M_\gamma d_VL_\beta , p)  \right) . \]
Recall the derivatives of $P^2= \left(\tr_\Sigma K\right)^2$ from Lemma \ref{lemma_con_check_1}.
\begin{align*}
d_M P^2 (X) &= 2 P \tr_\Sigma \nabla^M_X K \\
d_V P^2 (X) &= -4 P K(X, \nu) \\
\nabla^M_Y d_V P^2 (X) &= -4 \left( \tr_\Sigma \nabla^M_Y K \right) K(X,\nu) - 4 P \nabla^M_YK(X,\nu)
\end{align*}
We will calculate the three terms of $W$ separately. Since we choose normal coordinates around $p$, we have $p=0$. The relevant integrals are presented in Appendix \ref{chapter_spherical_integrals}. For better readability we drop the subscript from $\tr_M$.
\begin{enumerate}
  \item[1)]
	\begin{align*}
	c^{(\alpha,\beta)}\left((d_M P^2)_\beta, p \right) &= \int_{S_1} 2 P \tr_{S_1} \nabla^M_\beta K x^\alpha x^\beta \, \d\mu \\
	&= 2 	\int_{S_1}   \left(\tr K -  K(\nu,\nu)\right) \left( \tr \nabla^M_\beta K -\nabla^M_\beta K(\nu,\nu)\right)    x^\alpha x^\beta \, \d\mu
		\end{align*}
	
	\begin{align*}	
	\int_{S_1} 2 \tr K \tr \nabla^M_\beta K  x^\alpha x^\beta \, \d\mu &= 2 \tr K \tr \nabla^M_\beta K\int_{S_1}x^\alpha x^\beta  \, \d\mu \\
	&= \frac{8 \pi}{3}  \tr K \tr \nabla^M_\alpha K \\
	\int_{S_1} - 2 K(\nu,\nu) \tr \nabla^M_\beta K  x^\alpha x^\beta \, \d\mu &= -2 K_{\gamma \delta} \tr \nabla^M_\beta K \int_{S_1}x^\alpha x^\beta x^\gamma x^\delta  \, \d\mu \\
	&= - \frac{8 \pi}{15} \left( \tr K \tr \nabla^M_\alpha K + 2 \langle K(e_\alpha, \cdot) , \tr \nabla^M_\cdot K \rangle  \right)	\\
	\int_{S_1} -2 \tr K \nabla^M_\beta K(\nu,\nu) x^\alpha x^\beta    \, \d\mu &= -2 \tr K \nabla^M_\beta K_{\gamma \delta} \int_{S_1}  x^\alpha x^\beta x^\gamma x^\delta \, \d\mu \\
	&= - \frac{8 \pi}{15}  \left( \tr K \tr \nabla_\alpha^M K + 2 \tr K \Div_E K(e_\alpha)  \right)
	\end{align*}
	\begin{align*}
	\int_{S_1}  2 K(\nu,\nu) & \nabla^M_\beta K(\nu,\nu)    x^\alpha x^\beta  \, \d\mu = 2 K_{\gamma \delta} \nabla^M_\beta K_{\mu \nu} \int_{S_1} x^\alpha x^\beta x^\gamma x^\delta x^\mu x^\nu \, \d\mu \\
	&= \frac{8 \pi}{105}  \left( \tr K \nabla_\alpha \tr K  + 2 \tr K \Div K(e_\alpha,\cdot) + 2 \langle K(e_\alpha,\cdot), \nabla \tr K \rangle  \right.
\\
& \hspace{0.3cm} + 4 \langle K(e_\alpha,\cdot) , \Div K \rangle + \left. 2 \langle K, \nabla_\alpha K \rangle + 4 \langle K , \nabla K(e_\alpha, \cdot) \rangle \right) \\
&=: \frac{8 \pi}{ 105} I_\alpha
	\end{align*}
	
	\item[2)] 
	
	\begin{align*}
	c^{\beta}\left((d_M [d_V P^2]_\alpha)_\beta , p \right) &= -4 \int_{S_1} \left( \left( \tr_{S_1} \nabla^M_{e_\beta} K \right) K(e_\alpha,\nu) + P \nabla^M_{e_\beta} K(e_\alpha,\nu) \right) x^\beta \, \d\mu \\
	&= -4 \int_{S_1} \left( \left( \tr \nabla^M_{e_\beta} K -\nabla^M_{e_\beta} K (\nu,\nu)  \right) K(e_\alpha,\nu) \right. \\
	& \hspace{1cm} \left. + \left(  \tr K - K(\nu,\nu) \right) \nabla^M_{e_\beta} K(e_\alpha,\nu) \right) x^\beta \, \d\mu \\	
	\end{align*}
	
	\begin{align*}
	\int_{S_1}  -4 \tr \nabla^M_\beta K K(e_\alpha, \nu) x^\beta \, \d\mu &= -4 \tr K \nabla^M_\beta K_{\beta, \gamma} 	\int_{S_1}x^\beta x^\gamma \, \d\mu \\
	&= -\frac{16 \pi}{3} \tr \nabla^M_\beta K K_{\alpha \beta} \\
	&= -\frac{16 \pi}{3} \langle \tr \nabla^M_{\cdot} K , K(e_\alpha, \cdot) \rangle \\
	\int_{S_1} 4 \nabla_\beta^M K(\nu,\nu) K(e_\alpha,\nu) x^\beta \, \d\mu &= 4  \nabla_\beta^M K_{\gamma \delta} K_{\alpha \epsilon} 	\int_{S_1}x^\beta x^\gamma x^\delta x^\epsilon \, \d\mu \\
	&= \frac{16 \pi}{15} \left(  2 \langle \Div_E K , K(e_\alpha, \cdot) \rangle + \langle \tr\nabla^M_\cdot K , K(e_\alpha, \cdot)\rangle   \right) 	\\
\int_{S_1} -4 \tr K \nabla^M_{\beta} K(e_\alpha, \nu) x^\beta\, \d\mu &= -4 \tr K \nabla^M_{\beta} K_{\alpha \gamma} 	\int_{S_1} x^\beta x^\gamma \, \d\mu \\
&= - \frac{16 \pi}{3 } \tr K \Div_EK (e_\alpha)\\
\int_{S_1} 4 K(\nu,\nu) \nabla^M_\beta K(e_\alpha, \nu) x^\beta  \, \d \mu &= 4 K_{\gamma \delta} \nabla^M_\beta K_{\alpha \epsilon} \int_{S_1} x^\beta x^\gamma x^\delta x^\epsilon \, \d\mu \\
&= \frac{16 \pi}{15} \left( \tr K \Div_E K(e_\alpha) + 2 \langle K , \nabla^M_\cdot K(e_\alpha, \cdot) \rangle \right)
	\end{align*}
	
\item[3)]	

\begin{align*}
c^{(\alpha, \beta, \gamma)} \left( (d_M[d_V P^2]_\beta)_\gamma , p \right) &= -4 	\int_{S_1} \left( (\tr \nabla^M_\gamma K- \nabla^M_\gamma K(\nu,\nu) ) K(e_\beta,\nu) \right. \\
& \hspace{1cm} \left. + (\tr K - K(\nu,\nu))\nabla^M_\gamma K(e_\beta,\nu) \right) x^\alpha x^\beta x^\gamma  \, \d\mu
\end{align*}

\begin{align*}
	\int_{S_1}  -4 \tr \nabla^M_\gamma K K(e_\beta,\nu) x^\alpha x^\beta x^\gamma \, \d\mu &=  -4 \tr \nabla^M_\gamma K K_{\beta \delta} \int_{S_1}x^\alpha x^\beta x^\gamma x^\delta  \, \d\mu \\
	&= -\frac{16 \pi}{15} \left( \tr K \tr\nabla^M_\alpha K  + 2 \langle \tr \nabla^M_\cdot K ,K(e_\alpha, \cdot)  \rangle  \right) \\
	\int_{S_1} 4 \nabla^M_\gamma K(\nu,\nu) K(e_\beta, \nu) x^\alpha x^\beta x^\gamma   \, \d\mu &= 2 \frac{8 \pi}{ 105} I_\alpha \\
	\int_{S_1} -4 \tr K \nabla^M_\gamma K(e_\beta, \nu) x^\alpha x^\beta x^\gamma \, \d\mu &= -4 \tr K \nabla^M_\gamma K_{\beta \delta} \int_{S_1} x^\alpha x^\beta x^\gamma x^\delta \, \d\mu \\
	&= -\frac{16 \pi}{15} \left(  \tr K \tr \nabla^M_\alpha K  + 2 \tr K \Div_E K (e_\alpha) \right) \\
	\int_{S_1} 4 K(\nu,\nu)  \nabla^M_\gamma K(e_\beta, \nu) x^\alpha x^\beta x^\gamma   \, \d\mu &= 2 \frac{8 \pi}{ 105} I_\alpha \\
\end{align*}
\end{enumerate}

Adding all these terms up yields
\begin{align*}
W^\alpha_p [P^2] &= - \frac{8}{5} \langle \tr \nabla^M_\cdot K ,K( e_\alpha , \cdot)  \rangle  - \frac{8}{5} \tr K \Div_E K( e_\alpha ) + 4 \tr K \tr \nabla_\alpha K \\
& \hspace{0.3cm}- \frac{16}{5} \langle K( e_\alpha ,\cdot) , \Div_E K \rangle - \frac{16}{5}  \langle K , \nabla^M_\cdot K( e_\alpha ,\cdot) K \rangle + \frac{4}{5}  I_\alpha \\
&= \frac{4}{5}\partial_\alpha \left(  3 \tr K ^2 + |K|^2 \right) .
\end{align*}

\end{proof}

\begin{corollary}
Using the the various integrals calculated in the proof above we find
\begin{align*}
W^\alpha [ \tr K^2] &= 6 \partial_\alpha \tr K^2\\
W^\alpha [K(\nu,\nu)^2] &= \frac{2}{5} \partial_\alpha \left( \tr K ^2 + 2 |K|^2  \right)\\
W^\alpha [\tr K K(\nu,\nu)] &= 2 \partial_\alpha \tr K ^2 .
\end{align*}
In particular, this allows us to determine functions $L(x,\nu)$ such that $\HH_L$ concentrates at points with $\partial_\alpha \left( \Sc + \tr K^2 - |K|^2 \right) = 0$. For instance
\begin{align*}
L &= -\frac{3}{4} P^2 + 2 K(\nu, \nu)^2 ,   \\
L &= -\frac{1}{4} \tr K^2 + \frac{5}{4} K(\nu,\nu)^2 .
\end{align*}
\end{corollary}

Moreover, we can determine the expansion of $\HH$ on small spheres.
\begin{corollary}
Let $\Sigma \subset M$ be a spherical surfaces, with $|\Sigma| = 4 \pi R^2 $.  
Suppose $\Sigma$ is contained in a normal coordinate neighborhood $B_r(p)$ as in Lemma \ref{lemma_good_normal_coord} and that  $\| \A \|^2_{L^2(\Sigma)} \leq C_1 r R^2 $ for a constant $C_1$. 
Let $\HH_L$ be of Hawking type for an $L\in C^1$, then there is a constant $r_0 >0 $ such that for all $ r \in (0, r_0) $ we have the expansion
\[ \left| \HH_L[\Sigma] - 4\pi + \frac{2 \pi }{3} R^2 \Sc_p - R^2 c(L,p) \right| \leq C R^3 . \]
Where $C$ depends only on $r_0$, $C_1$, $L$  and $C_B$. 
Additionally, we calculate 
\begin{align*}
c(\tr K^2,p ) &=  4 \pi  \tr K_p^2 \\
c(K(\nu,\nu)^2,p) &= \frac{4\pi}{15} \left( \tr K_p^2 + 2 |K_p|^2 \right) \\
c(\tr K K(\nu,\nu),p) &= \frac{4 \pi}{ 3} \tr K_p^2.
 \end{align*}
Hence we obtain
\[ \left| \HH[\Sigma] - 4\pi + \frac{2 \pi}{3} R^2 \Sc_p + \frac{ 2 \pi }{15} R^2 (3 \tr K_p^2 + |K_p|^2)  \right| \leq C R^3  \]
for the functional corresponding to the Hawking energy.
\end{corollary}

\begin{proof}

From \cite[Theorem 5.1]{LMI} we get the expansion
\[ \left| \WW[\Sigma] - 4\pi + \frac{2 \pi }{3} R^2 \Sc_p \right| \leq C R^3 . \]
Using the coordinates of Lemma \ref{lemma_good_normal_coord} and Lemma \ref{cor_sigma_to_s} we have $\LL[\Sigma] = R_E^2 c(L,p) + O(R^3)$.
This implies the expansion since $R$ and $R_E$ are all comparable. The explicit calculation of $c$ is straight forward, using the results of Appendix \ref{chapter_spherical_integrals}. 
\end{proof}

\begin{corollary}
Let $L = \alpha \tr K^2 + \beta K(\nu,\nu)^2 + \gamma \tr K K(\nu,\nu) $ for $\alpha,\, \beta,\, \gamma \in \R$, then we have 
\[ W_p[L] = \nabla^M \frac{3}{2 \pi} c(L,p), \]
i.e. a concentration point of this kind of Hawking type functional is a critical point of its second order expansion.
\end{corollary}

\newpage


\begin{appendices}

\section{Bounded Geometry and Small Surfaces}
\label{chapter_bounded_geo_small_surf}

In this section we will briefly introduce bounded geometry, as presented in \cite{LMI} and \cite{LMII}. For a more comprehensive treatment see for instance \cite[Chapter 2]{Eldering_NHIM}.

\begin{definition}
Let $(M,g)$ be a complete Riemannian manifold with injectivity radius $r_{inj}(M,g,p)$ at $p\in M$ and Riemannian curvature tensor $\Rm$. We say $M$ has $C_B$-bounded geometry if there exists a constant $C_B > 0$ such that for each $p \in M$ we have 
\[ r_{inj}(M,g,p) \geq C_B^{-1} \]
and 
\[ |\Rm| + |\nabla \Rm| \leq C_B . \] 
\end{definition}

We may combine the well known results on normal coordinates with the uniform bound on the injectivity radius to obtain the following lemma.
\begin{lemma}[cf. {\cite[Section 2.1]{LMI}}]
Let $(M,g)$ be a manifold of $C_B$ bounded geometry, let $B_r(y)$ be the Euclidean ball at $y\in T_p M $ of radius $r$ and $\BB_r(p)$ the geodesic ball at $p\in M$ with radius $r$. There exist constants $h_0$ and $r_0$, depending only on $C_B$, such that in normal coordinates $\exp_p: B_{r_0}(0) \to \BB_{r_0}(p)$ the metric satisfies 
\[ g= g_E + h ,\]
where $g_E$ is the Euclidean metric and $h$ obeys
\[ \sup_{B_{r_0}(0)} \left(|x|^2 |h| + |x||\nabla^E h| + \left|\left(\nabla^E\right)^2 h \right| \right) \leq h_0 .\]
Here $x$ is the position vector field in $B_r(0)$, $|\cdot|$ is the Euclidean norm and $\nabla^E$ is the Euclidean connection.
\end{lemma}

Next we consider small surfaces $(\Sigma, \gamma)$ that are isometrically immersed in a three dimensional, $C_B$ bounded manifold $(M,g)$. 
That is we deal with closed surfaces contained in geodesic balls $\Sigma \subset \BB_{r_0}(y)$ for some point $y \in M$ and $r_0 \leq \min (r_{inj}, 1)$. 
With our previous result in mind, we regard them as immersed in $B_{r_0}(0)$ equipped with the metric $g= g_E + h$ as above. We fix this setting for now, unless stated otherwise.
Additionally, we will denote all geometric quantities computed with respect to the Euclidean metric by a $E$ index.

\begin{lemma}[see {\cite[Lemma 2.1]{LMI}}]
\label{lemma_small_surf_euc}
There exist a constant $C$, depending only on $r_0$ and $h_0$, such that for all surfaces $\Sigma \subset B_r(0)=B_r$, $r \leq r_0$, we have
\begin{align*}
|\gamma - \gamma_E|_E &\leq C |x|_E^2 \\
| \sqrt{|\det \gamma |} - \sqrt{|\det \gamma_E |} | & \leq C \sqrt{ \det \gamma_E} |x|^2_E \\ 
| \sqrt{|\det \gamma |} - \sqrt{|\det \gamma_E |} | & \leq C \sqrt{ \det \gamma} |x|^2_E \\ 
|\nu - \nu_E| & \leq C |x|_E^2 \\
|A - A_E|_E & \leq C (|x|_E + |x|_E^2 |A_E|_E)
\end{align*}
\end{lemma}

\begin{definition}
\label{def_radius}
We define the radius $R$ of $\Sigma$ with respect to $\gamma$ by the relation $|\Sigma|=:4\pi R^2$. Analogously, the corresponding Euclidean radius $R_E$ is given by $|\Sigma|_E =:4 \pi R_E^2$, where $|\Sigma|_E := \int_\Sigma \d \mu_{\gamma_E}$.
\end{definition}

\begin{corollary}[cf. {\cite[Lemma 2.5]{LMI}}]
\label{cor_small_surf_euc}
In the setting of Lemma \ref{lemma_small_surf_euc} we have that the following estimates.
\begin{enumerate}
	\item 
	\begin{align*}
	\left| |\Sigma| - |\Sigma|_E  \right| &\leq  C r^2 |\Sigma| \\
		\left| |\Sigma| - |\Sigma|_E  \right| &\leq  C r^2 |\Sigma|_E \\
		|R - R_E| &\leq C r^2 R  \\
		|R - R_E| &\leq C r^2  R_E
	\end{align*}
	In particular, the areas $|\Sigma|$ and $|\Sigma|_E$ are comparable, as are the corresponding radii $R$  and $R_E$.
		
	\item 
\begin{align*}
\left| \WW[\Sigma,g] - \WW[\Sigma, g_E] \right| &\leq C(C_B) r^2 \left(|\Sigma| + \WW[\Sigma,\gamma] + r^2 \|\A \|^2_{L^2(\Sigma,\gamma)} \right) \\
\| \A^E \|^2_{L^2(\Sigma,\gamma_E)} &\leq C(C_B) \| \A \|^2_{L^2(\Sigma,\gamma)} + C r^4 \WW[\Sigma,g]
\end{align*}

\end{enumerate}

\end{corollary}

\begin{lemma}[see {\cite[Lemma 2.2]{LMI}}]
\label{lemma_area}
There exists $0< r_1 \leq r_0$ and a purely numerical constant $C$ such that for all $\Sigma \subset B_r$, $r \leq r_1$, we have 
\[ |\Sigma| \leq C r^2 \WW[\Sigma]. \]
\end{lemma}
\begin{lemma}[see {\cite[Lemma 2.5]{LMII}}]
\label{lemma_diameter}
There exists a constant $C$, depending only on $C_B$, such that all connected surfaces $\Sigma \subset M$ obey
\[ \diam_M(\Sigma) \leq C \left(|\Sigma|^{1/2} \WW[\Sigma]^{1/2} + |\Sigma| \right) . \]
\end{lemma}
Clearly, the previous two lemmas also hold for stratified surfaces if we apply them to every component.

In Section  \ref{section_surface_concentration} it is necessary to approximate a surfaces by a sphere, hence we state a scaled version of the results of De Lellis and Müller on that topic together with an estimate on the normal vectors.

\begin{theorem}[cf {\cite[Theorem 1.1]{DeLellis_Mueller_rigidity}} and {\cite[Theorem 1.2]{DeLellis_Mueller_umbilical}}]
\label{thm_sigma_to_s}
Let $\Sigma \subset \R^3$ be a surfaces with induced metric $\gamma_E$ and $\| \A_E \|_{L^2(\Sigma, \gamma_E)}^2 < 8 \pi$ and consider its Euclidean radius $R_E$ as well as its Euclidean center of gravity $a_E:=|\Sigma|_E^{-1} \int_\Sigma x \, d \mu_{\gamma_E}$. Then there exists a universal constant $C$ and a conformal map $\psi : S := S_{R_E}^2(a_E) \to \Sigma$ with the following properties. Let $\sigma$ be the round metric on $S$, $\nu_S$ its unit normal vector field and  let $\alpha$ be the conformal factor of $\psi$, i.e. $\psi^*\gamma_E = \alpha^2 \sigma$. Then the following estimates hold.
\begin{align*}
\left\| A_E - R_E \id \right\|_{L^2(\Sigma,\gamma_E)} &\leq C \|\A_E\|_{L^2(\Sigma,\gamma_E)}\\
\left\| H_E - \frac{2}{R_E} \right\|_{L^2(\Sigma,\gamma_E)} &\leq C \|\A_E\|_{L^2(\Sigma,\gamma_E)} \\
\| \psi -(a_E + \id_S) \|_{L^2(S)} & \leq C R_E^2 \|\A_E\|_{L^2(\Sigma,\gamma_E)} \\
\| \psi -(a_E + \id_S) \|_{L^\infty(S)} & \leq C R_E \|\A_E\|_{L^2(\Sigma,\gamma_E)} \\
\| \d \psi - \id_{TS} \|_{L^2(S)} & \leq C R_E \|\A_E\|_{L^2(\Sigma,\gamma_E)} \\
\| \alpha - 1 \|_{L^\infty(S)} & \leq C \|\A_E\|_{L^2(\Sigma,\gamma_E)} 
\end{align*}
\end{theorem}

\begin{corollary}
\label{cor_normal_sigma_s}
Assume additionally that $\|\A_E\|_{L^2(\Sigma,\gamma_E)}$ is so small that  $\| \alpha - 1 \|_{L^\infty(S)}  \leq 1/2$. Then there is a universal constant $C$ such that
\[ \| \nu_S - \nu_E \circ \psi \|_{L^2(S)}  \leq  C R_E \|\A_E\|_{L^2(\Sigma,\gamma_E)} . \]
\end{corollary}

We may combine the results of Lemma \ref{lemma_small_surf_euc}, Theorem \ref{thm_sigma_to_s} and the previous corollary in order to approximate a small surface $(\Sigma, \gamma)$ with $(S, \sigma)$.

\begin{corollary}
\label{cor_sigma_to_s}
Let $\Sigma \subset B_{r_0} $ be a small surface, and assume that $\|\A_E\|_{L(\Sigma,\gamma_E)}  \leq 8\pi$ is small enough that corollary \ref{cor_normal_sigma_s} holds. Assume further that $H^{-1}$ and $H_E^{-1}$ are uniformly bounded. Then we have the following estimates for a constant $C$ dependent only on $C_B$.
\begin{align*}
| \d \, \mu_\gamma - \d \, \mu_\sigma  | & \leq C \left(r^2 + \| \A_E\|_{L^2(\Sigma,\gamma_E)}\right) \\
\| \nu \circ \psi - \nu_S \|_{L^2(S,\sigma)} & \leq C \left( r^2 R_E + R_E \|\A_E\|_{L^2(\Sigma,\gamma_E)}  \right) \\
\| H^{-1} - R_E/2 \|_{L^2(\Sigma,\gamma_E)} & \leq C\left( \sup_\Sigma \frac{r }{H H_E} \left(1 + r^2 \| A_E \|^2_{L^2(\Sigma,\gamma_E)}\right) +\sup_\Sigma \frac{R_E}{H_E} \| \A_E \|_{L^2(\Sigma, \gamma_E)}\right)
\end{align*}
Moreover,  we may transport any bounded Lipschitz function $F: \R^3 \times \R^3 \to \R$ from $\Sigma$ to $S$.
 \[\left|\int_\Sigma F(y,\nu) \d\,\mu_\gamma - \int_S F(x,\nu_S) \d\,\mu_\sigma \right| \leq C_1 R_E^2 \left(r^2 + \| \A_E\|_{L^2(\Sigma,\gamma_E)}\right) \]
Here $C_1$ is a constant that depends on $C_B$ and $F$.
\end{corollary}

It is possible to choose normal coordinates of $(M,g)$ well suited for a given closed surface $\Sigma$.

\begin{lemma}[see {\cite[Lemma 3.1]{Metzger_refined_position}}]
\label{lemma_good_normal_coord}
Let $\Sigma \subset M$ be a surface with extrinsic diameter $d$ such that $2d \leq \inj(M,g)$. Then there exists a point $p_0\in M$ with $\dist(p_0,\Sigma ) \leq d$ and such that in normal coordinates $\psi$ centered at $p_0$ we have that
\[ a = \frac{1}{|\Sigma|} \int_{\psi(\Sigma)} y \,\d\mu_g = 0\]
and
\[  |a_E|_E = \frac{1}{|\Sigma|_E}\left| \int_{\psi(\Sigma)} y \,\d\mu_E \right|_E \leq  C(C_B) d^3 . \]
Additionally, if $\Sigma$ obeys $\| \A^E \|^2_{L^2(\Sigma, \gamma_E)} \leq 8\pi$ then we have
\[
\max_{x \in \Sigma} |x|_E \leq C(C_B) R_E .
\]
\end{lemma}


\section{Useful Integrals over the Sphere}
\label{chapter_spherical_integrals}

Consider the 2-sphere $S:=S^2_1(0)$ and let $\{ x^{\alpha} \}$ be the Euclidean coordinates on in $\R^3$. We calculate integrals of the form
 \[ \int_S x^{\alpha_1}... \, x^{\alpha_n} \, d \mu^S,  \]
where $n$ is a natural number up to $6$. If $n$ is odd, the integral always vanishes. In the other cases a lengthy but straight forward calculation reveals the following.
\begin{align*}
\int_S x^\alpha x^\beta \, d\mu &= \frac{4 \pi}{3} \delta^{\alpha \beta} \\
\int_S x^\alpha x^\beta x^\gamma x^\delta \, d\mu &= \frac{4 \pi}{15} \left( \delta^{\alpha \beta}\delta^{\gamma \delta} + \delta^{\alpha \gamma }\delta^{\beta \delta} +\delta^{\alpha \delta}  \delta^{\beta \gamma} \right)\\
\int_S x^\alpha x^\beta x^\gamma x^\delta x^\epsilon x^\rho \, d\mu &= \frac{4 \pi}{105} \left(\delta^{\alpha \beta}\delta^{\gamma \delta} \delta^{\epsilon \rho} +  \delta^{\alpha \beta}\delta^{\gamma \epsilon} \delta^{\delta \rho} +  \delta^{\alpha \beta}\delta^{\gamma \rho} \delta^{\epsilon \delta}  \right. \\
& \hspace{0.3cm} + \delta^{\alpha \gamma}\delta^{\beta \delta} \delta^{\epsilon \rho} + \delta^{\alpha \gamma}\delta^{\beta \epsilon} \delta^{\delta \rho} + \delta^{\alpha \gamma}\delta^{\beta \rho} \delta^{\epsilon \delta}\\
& \hspace{0.3cm} + \delta^{\alpha \delta}\delta^{\gamma \beta} \delta^{\epsilon \rho} +\delta^{\alpha \delta}\delta^{\gamma \epsilon} \delta^{\beta \rho} +\delta^{\alpha \delta}\delta^{\gamma \rho} \delta^{\epsilon \beta} \\
& \hspace{0.3cm} + \delta^{\alpha \epsilon}\delta^{\gamma \delta} \delta^{\beta \rho} + \delta^{\alpha \epsilon}\delta^{\gamma \beta} \delta^{\delta \rho} +\delta^{\alpha \epsilon}\delta^{\gamma \rho} \delta^{\beta \delta} \\
& \hspace{0.3cm} \left. + \delta^{\alpha \rho}\delta^{\gamma \delta} \delta^{\beta \epsilon} + \delta^{\alpha \rho}\delta^{\gamma \beta} \delta^{\delta \epsilon} + \delta^{\alpha \rho}\delta^{\gamma \epsilon} \delta^{\beta \delta} \right)
\end{align*}

\end{appendices}

\bibliographystyle{hplain} 
\bibliography{thesis_bib} 

\end{document}